\documentclass[12pt]{article}
\usepackage{latexsym}
\newtheorem{theorem}{Theorem}
\newcommand{\eqdef}{\, =\kern -11pt\raise 6pt\hbox{{\tiny\textrm{def}}}\,}
\begin{document}
\title{Identically Distributed Pairs of Partition Statistics}
\author{Herbert S. Wilf\\University of Pennsylvania, Philadelphia, PA 19104-6395}
\date{}
\maketitle
\begin{abstract}
We show that many theorems which assert that two kinds of partitions of the same integer $n$ are equinumerous are actually special cases of a much stronger form of equality. We show that in fact there correspond partition statistics $X$ and $Y$ that have identical distribution functions. The method is an extension of the principle of sieve-equivalence, and it yields simple criteria under which we can infer this identity of distribution functions.
\end{abstract}
Let ${\cal P}(n)$ be the set of partitions of the integer $n$, $p(n)=|{\cal P}(n)|$, and let ${\cal P}=\cup_{n\ge 1}{\cal P}(n)$ be the set of all partitions of all positive integers. A \textit{partition statistic} $X$ is a nonnegative-integer-valued function defined on ${\cal P}$. Two partition statistics $X$ and $Y$ are \textit{identically distributed} if 
\[\forall n\ge 1,\forall j\ge 0: \mathrm{Prob}_n(X=j)=\mathrm{Prob}_n(Y=j),\]
where
\[\mathrm{Prob}_n(X=j)\eqdef\ \frac{|\{\pi\in{\cal P}(n):X(\pi)=j\}|}{p(n)}.\]

There are a number of classical theorems of the form ``The number of partitions of $n$ that have no \dots is the same as the number that have no \dots.'' The purpose of this note is to remark that frequently in such cases there is an underlying pair of identically distributed partition statistics. For example, Euler's famous theorem that the number of partitions of $n$ with distinct parts is the same as the number with odd parts can be strengthened to the following.
\begin{theorem}
\label{th:A}
The partition statistics $X(\pi)=$ number of even part sizes that occur in $\pi$, and $Y(\pi)=$ number of repeated part sizes that occur in $\pi$, are identically distributed.
\end{theorem}
That is, for every $n$ and $j$, the number of partitions of $n$ in which exactly $j$ different sizes of parts are repeated is the same as the number in which exactly $j$ different sizes of parts are even, and Euler's original theorem is the case $j=0$. Thus a great deal more is true than is stated in the original theorem.

Let's prove Theorem \ref{th:A}, for then the general idea will be clear. Imagine that a sieve-method computation is carried out in which the objects are the partitions of $n$. Suppose the $r$th property of a partition is that it contains the pair of parts $\{r,r\}$. Now imagine that another sieve-method calculation is carried out on the same objects, in which now the $r$th property of a partition is that it contains the single part $2r$. We claim that the output from these two sieve calculations will be the same, i.e., that the two computations are \textit{sieve-equivalent} \cite{W}.

Indeed, we claim that the inputs are the same, so the outputs must be the same too. The input to the first sieve computation is, for every $S$, the number of partitions that have at least the multiset $S$ of repeated parts $\{r_1,r_1,r_2,r_2,\dots,r_k,r_k\}$, and that number is $p(n-2r_1-2r_2-\dots -2r_k)$. But this is the same as the input to the second sieve calculation. Hence both sieves give the same outputs. Among the outputs of a sieve calculation are the numbers of objects that have exactly $j$ properties, for each $j$. Hence these numbers are the same, and Theorem \ref{th:A} is proved. $\Box$

We now give, following Cohen \cite{C} and Remmel \cite{R}, two general theorems that imply the identical distribution of a pair of partition statistics. The first theorem is easier to work with, but the second is more general.
\begin{theorem}
\label{th:B}
Let $\{F_1,F_2,\dots \}$, and $\{G_1,G_2,\dots \}$ be two lists of multisets of positive integers. Suppose that the $F$'s are pairwise disjoint, the $G$'s are pairwise disjoint, and that for every $i=1,2,\dots$ we have $\sum_{k\in F_i}k=\sum_{k\in G_i}k$. Define the following two partition statistics:
\[X(\pi)=|\{i:F_i\in\pi\}|;\quad Y(\pi)=|\{i:G_i\in\pi\}|.\]
Then $X$ and $Y$ are identically distributed.
\end{theorem} 

This is the so-called ``disjoint-multiset'' case. Here are some other pairs of identically distributed partition statistics that follow from the theorem above.
\begin{enumerate}
\item $X=$ the number of part sizes that are perfect squares; $Y=$ the number of part sizes $i$ whose multiplicity is $\ge i$.
\item $X=$ the number of part sizes that are $\equiv$ 2,3,4 mod 6; $Y=$ the number of part sizes that are either a multiple of 3 or else repeated and not a multiple of 3.
\item Fix an integer $d>1$. Let $X=$ the number of part sizes that are multiples of $d$; $Y=$ the number of part sizes whose multiplicity is $\ge d$.
\item Let $M_1, M_2$ be two sets of positive integers. Define $2M_1=\{j:(j/2)\in M_1\}$. Suppose that $2M_1\subseteq M_1$ and $M_2=M_1-2M_1$. Then define the statistic $X=$ the number of part sizes that are not in $M_2$; $Y=$ the number of part sizes $i$ such that either $i\notin M_1$ or $i\in M_1$ and is repeated.
\end{enumerate}

These last three examples generalize theorems of Schur \cite{S}, Glaisher \cite{G}, and Andrews \cite{A}, respectively. They were originally used by Remmel \cite{R} as examples of his elegant disjoint multisets machinery for producing partition bijections. 

The hypotheses of Theorem \ref{th:B}, which require pairwise disjointness, can be weakened to sieve-equivalence. We state this as follows.
\begin{theorem}
\label{th:C}
Let $\{F_1,F_2,\dots \}$, and $\{G_1,G_2,\dots \}$ be two lists of multisets of positive integers. For every finite set $S$, put ${\cal F}_S=\cup_{j\in S}F_j$, in the sense of a multiset union, and similarly for ${\cal G}_S$. Suppose that for every $S$ we have $\sum_{j\in {\cal F}_S}j=\sum_{j\in {\cal G}_S}j$. Then, with $X,Y$ defined as in Theorem \ref{th:B}, $X$ and $Y$ are identically distributed.
\end{theorem}

As an application of Theorem \ref{th:C} which does not follow from Theorem \ref{th:B} we can use the pair of lists of multisets, again taken from \cite{R},
\[\{\{2,4\},\{4,6\},\{6,8\},\dots\}\qquad \{\{1,1,2,2\},\{2,2,3,3\},\{3,3,4,4\},\dots\}\]
to deduce that the following partition statistics are identically distrbuted: $X=$ the number of consecutive even part sizes, and $Y=$ the number of consecutive repeated part sizes. This is a considerable sharpening of Corollary 2.10 of \cite{R}, and many other examples are easy to construct.

\end{document}